\documentclass[11pt]{amsart}
\usepackage{amssymb, paralist, xspace, graphicx, url, amscd, euscript, mathrsfs,stmaryrd,epic,eepic,color}
\usepackage[all]{xy}
\SelectTips{cm}{}

\numberwithin{equation}{section}

\numberwithin{subsection}{section}

\allowdisplaybreaks[1]

\newtheorem*{namedtheorem}{\theoremname}
\newcommand{\theoremname}{testing}

\newtheorem{theorem}[subsection]{Theorem}

\newtheorem{proposition-definition}[subsection]
{Proposition-Definition}

\theoremstyle{definition}
\newtheorem{definition}[subsection]{Definition}

\theoremstyle{remark}

\newcommand\cI{\mathcal{I}}

\newcommand\cO{\mathcal{O}}

\newcommand\GG{\mathbb{G}}

\newcommand\bP{\mathbf{P}}

\newcommand\arr{\ifinner\to\else\longrightarrow\fi}

\newcommand\dar{\downarrow}
\newcommand\rar{\rightarrow}
\newcommand\das{\dashrightarrow}

\def\displaytimes_#1{\mathrel{\mathop{\times}\limits_{#1}}}

\def\displayotimes_#1{\mathrel{\mathop{\bigotimes}\limits_{#1}}}

\newcommand\Aut{\operatorname{Aut}}

\newdir{ >}{{}*!/-5pt/@{>}}

\newcommand\doublelong[2]{\mathbin{\xymatrix{{}\ar@<3pt>[r]^{#1}
\ar@<-3pt>[r]_{#2}&}}}

\newlength{\ignora}

\renewcommand{\setminus}{\smallsetminus}
\newcommand{\setmin}{\smallsetminus}

\theoremstyle{plain}

\newtheorem*{thm*}{Theorem}

\newtheorem{lem}[subsection]{Lemma}

\theoremstyle{definition}

\numberwithin{equation}{subsection}

\begin{document}

\title{Weak toroidalization over non-closed fields}

\author[Abramovich]{Dan Abramovich}

\author[Denef]{Jan Denef}

\author[Karu]{Kalle Karu}

\address[Abramovich]{Department of Mathematics\\
Brown University\\
Box 1917\\
Providence, RI 02912\\
U.S.A.}
\email{abrmovic@math.brown.edu}

\address[Denef]{University of Leuven, Department of Mathematics, Celestijnenlaan 200 B,
B-3001 Leuven, Belgium.}
\email{jan.denef@wis.kuleuven.be  }

\address[Karu]{Dept. of Mathematics
University of British Columbia
1984 Mathematics Rd
Vancouver, B.C. V6T 1Z2}
\email{karu@math.ubc.ca}

\thanks{Abramovich support in part by NSF grants, DMS-0901278. Karu was supported by NSERC grant RGPIN 283301}

\date{March 26, 2013}

\begin{abstract}
We prove that any dominant morphism of algebraic varieties over a field $k$ of characteristic zero can be transformed into a toroidal (hence monomial) morphism by projective birational modifications of source and target. This was previously proved by the first and third author when $k$ is algebraically closed. Moreover we show that certain additional requirements can be satisfied.
\end{abstract}

\maketitle

\setcounter{tocdepth}{1}


\section{Introduction}

The purpose of the present paper is to extend the toroidalization theorem obtained in \cite{Abramovich-Karu} by the first and third author, to the case of varieties over fields of characteristic zero that are not necessarily algebraically closed, and to slightly sharpen it by requiring that the modification of the source variety is an isomorphism outside the toroidal divisor. This extension is essential for a recent paper \cite{Denef-AxKochen} of the second author. Our result is the following theorem.
\begin{theorem}\label{th-toroidal-reduction}
Let $k$ be a field of characteristic zero. Let $f: X\to B$ be a dominant morphism
of varieties over $k$, and let $Z\subset X$ be a proper closed subset.
Then there exists a commutative diagram
$$\begin{array}{lclcl} U_{X' }& \subset & X' &\stackrel{m_X}{\to} &X \\
                     \dar & & \dar f' & & \dar f\\ U_{B'} & \subset & B'
                     & \stackrel{m_B}{\to} &B \end{array}
$$
where $X'$ and $B'$ are nonsingular quasi-projective varieties over $k$, and $m_B$ and $m_X$ are projective birational  morphisms, such that
\begin{enumerate}
\item the inclusions on the left are strict toroidal embeddings;
\item $f'$ is a toroidal quasi-projective morphism;
\item let $Z' = m_X^{-1}(Z)$, then $Z'$ is a strict normal crossings divisor, and $Z'\subset X'\setmin U_{X'}$;
\item the restriction of the morphism $m_X$ to  $U_{X' }$ is an open embedding.
\end{enumerate}
\end{theorem}
\noindent Note that when $f$ is proper, so is $f'$, hence $f'$ becomes projective.

Toroidal embeddings and toroidal morphisms are defined in Section~\ref{sec-notation} below. Because both varieties $X'$ and $B'$ are nonsingular, the embeddings $U_{X'}  \subset X'$, $U_{B'}  \subset  B'$ and the morphism $f'$ of the theorem can be described as follows. The requirement that the embeddings are strict toroidal is equivalent to the statement that $X'\setmin U_{X'}$,  $B'\setmin U_{B'}$ are strict normal crossings divisors.  The requirement that $f'$ is toroidal is equivalent to the following: after base change to an algebraic closure $\bar{k}$ of $k$, for each closed point $x\in X'_{\bar{k}}$, $b=f'(x)\in B'_{\bar{k}}$, there exist uniformizing parameters $x_1,\ldots,x_n$ for $\hat{\cO}_{ X'_{\bar{k}, x}}$ and $b_1,\ldots,b_m$ for $\hat{\cO}_{ B'_{\bar{k}, b}}$, such that
\begin{enumerate}
\item Locally at $x$, the product $x_1\cdots x_n$ defines the divisor $X'_{\bar{k}} \setmin U_{X', \bar{k}}$.
\item Locally at $b$, the product $b_1\cdots b_m$ defines the divisor $B'_{\bar{k}} \setmin U_{B', \bar{k}}$.
\item The morphism $f'$ gives $b_i$ as monomials in $x_j$.
\end{enumerate}

Here we say that elements $z_1,\ldots,z_n$ of a local ring $A$ containing its residue field are uniformizing parameters for $A$ if these elements minus their residues form a system of regular parameters for $A$.

Theorem \ref{th-toroidal-reduction} is related to Cutkosky's local monomialization theorem \cite{Cutkosky-LocMon}, which states that any dominant morphism of complete varieties over $k$, can be transformed by birational transformations (of source and target), which are products of monoidal transforms, into a monomial morphism of nonsingular quasi-complete (not necessarily separated) integral schemes of finite type over $k$.

Our proof of Theorem \ref{th-toroidal-reduction} does not depend on Cutkosky's Theorem, but follows very closely the proof of Theorem 2.1 in \cite{Abramovich-Karu}, which in turn uses methods and results from \cite{AJ, deJong-families}. The present paper would be much shorter if we would only elaborate what has to be upgraded in \cite{Abramovich-Karu} to obtain our theorem. But this would result in a very terse proof whose correctness is difficult to verify. Instead, we tried to make the present paper quite self-contained, repeating arguments from \cite{Abramovich-Karu, AJ}, except that we use the semistable reduction theorem of de Jong \cite{deJong-families}.

A recent closely related manuscript of Illusie and Temkin \cite{Illusie-Temkin} goes quite a bit farther than our Theorem \ref{th-toroidal-reduction} in that it proves Gabber's version of de Jong's alteration theorem in great generality; it thus requires much more technique. Their statement in characteristic zero is \cite[Theorem 3.9]{Illusie-Temkin}, which is similar to our Theorem \ref{th-toroidal-reduction} but is stated in the generality of quasi-excellent schemes of characteristic zero. A key difference from our method is that they do not use the torific ideal.

Note that Theorem \ref{th-toroidal-reduction} is much weaker than the strong toroidalization conjecture (see section 6.2 of \cite{AKMW}), which has been proved by Cutkosky \cite{Cutkosky-Tor3folds} for dominant morphisms of varieties of dimension $\leq 3$ over an algebraically closed field of characteristic zero.

\section{Notations and definitions}\label{sec-notation}

We work over a field $k$ of characteristic zero.  A \emph{variety (defined) over $k$} is an integral separated scheme of finite type over $k$. If $X$ is a variety defined over $k$, we let $X_{\bar{k}}$ be its base extension to an algebraic closure $\bar{k}$ of $k$. Note that the scheme $X_{\bar{k}}$ might not be a variety.

A {\em modification} is a proper birational morphism of varieties. An {\em alteration} is a proper, surjective, generically finite morphism of varieties. An alteration $Y\to X$ is called a {\em Galois alteration} with Galois group $G$ if the function field extension  $K(Y)/K(X)$ is Galois with Galois group $G$, and if the action of $G$ on $K(Y)$ is induced  by an action on $Y$ keeping the morphism $Y\to X$ invariant.

\subsection{Divisors}

Let $X$ be a smooth variety defined over $k$ (or more generally a smooth scheme over $k$), and $D\subset X$ a divisor. We say that $D$ is a {\em strict normal crossings divisor} if for every point $x\in X$ there exists a regular system of parameters $z_1,\ldots,z_n$ at $x$, such that every irreducible component of $D$ containing $x$ has local equation $z_i=0$ for some $i$. A divisor is a normal crossings divisor if it becomes a strict normal crossings divisor on some \'etale cover of $X$. The condition of being a (strict) normal crossings divisor is stable under base extension to algebraic closure: if $D\subset X$ is a (strict) normal crossings divisor, so is $D_{\bar{k}} \subset X_ {\bar{k}}$.

Let a finite group $G$ act on a (not necessarily smooth) variety $X$ over $k$ (or more generally a scheme of finite type over $k$), mapping a divisor $D\subset X$ into $D$. We say that $D$ is {\em $G$-strict} if the union of translates of each irreducible component of $D$ is normal. In the case where $G$ is the trivial group $1$, we say that $D$ is a strict divisor. Thus, a strict normal crossings divisor is both strict and normal crossings divisor.

\subsection{Toroidal embeddings}\label{sec-toroidal-embeddings}

We refer to \cite{KKMS, Fulton} for details about toric varieties. If $V$ is a toric variety, we denote by $T_V\subset V$ the big algebraic torus in $V$. Toric morphisms are always assumed to be dominant.

An open embedding of varieties $U\subset X$ defined over $\bar{k}$ (or more generally of schemes of finite type over $\bar{k}$)  is called a {\em toroidal embedding} if for every closed point $x\in X$ there exists a toric variety $V$ over $\bar{k}$, a closed point $v\in V$, and an isomorphism of complete local $\bar{k}$-algebras:
\[ \hat{\cO}_{X,x}  \cong  \hat{\cO}_{V,v},\]
such that the completion of the ideal of $X\setmin U$ maps isomorphically to the completion of the ideal of $V\setmin T_V$. The pair $(V,v)$, together with the isomorphism, is called a local model at $x\in X$. A toroidal embedding $U\subset X$ over $\overline{k}$ is called {\em strict} if $D=X\setmin U$ is a strict divisor.

An open embedding $U\subset X$ defined over $k$ is called a toroidal embedding if the base extension $U_{\bar{k}}\subset X_{\bar{k}}$ is a toroidal embedding. The toroidal embedding is called strict if the divisor $D=X\setmin U$ is strict. Note that if the toroidal embedding $U\subset X$ is strict, then the toroidal embedding $U_{\bar{k}}\subset X_{\bar{k}}$ is also strict, but the converse may not hold.

Let $U_X\subset X$ and $U_B\subset B$ be two toroidal embeddings defined over $\bar{k}$, and let $f:X\to B$ be a dominant morphism mapping $U_X$ to $U_B$. (We write such a morphism as $f:(U_X\subset X) \to (U_B\subset B)$.) Then $f$ is called a {\em toroidal morphism} if for every closed point $x\in X$ there exist local models $(V,v)$ at $x\in X$ and $(W,w)$ at $f(x)\in B$, and a toric morphism $g: V\to W$, with $f(v)=w$, such that the following diagram commutes:
\[
\begin{CD}
\hat{\cO}_{X,x} @>{\cong}>> \hat{\cO}_{V,v}\\
@A{\hat{f}^{\#}}AA @AA{\hat{g}^{\#}}A\\
\hat{\cO}_{B,f(x)} @>{\cong}>> \hat{\cO}_{W,w}
\end{CD}
\]
Here $\hat{f}^{\#}$ and $\hat{g}^{\#}$ are the ring homomorphisms coming from $f$ and $g$. In this situation we say that the local models $(V,v)$, $(W,w)$ are \emph{compatible with} $f$.

A morphism $f:(U_X\subset X) \to (U_B\subset B)$ between toroidal embeddings defined over the field $k$ is called a toroidal morphism if its base extension to $\bar{k}$ is a toroidal morphism.

The composition of two toroidal morphisms is again toroidal \cite{Abramovich-Karu}.

\subsubsection{Remark} \label{remark-rationality}
When $f:(U_X\subset X) \to (U_B\subset B)$ is a toroidal morphism of strict toroidal embeddings defined over $k$, and if $x$ is a closed point of $X$ with residue field $k(x)$ equal to $k$, or with $k(x)$ algebraically closed, then there exist local models at $x\in X$ and at $f(x)\in B$, which are compatible with $f$ and
\textquotedblleft defined\textquotedblright $\,$ over $k$. We refer to \cite{Denef-ToroidalMorphisms} for a precise statement of this assertion. This can be proved by adapting the argument in section 3.13 of \cite{Kato-LogStruct}. A detailed self-contained proof and some related results can be found in \cite{Denef-ToroidalMorphisms}. Moreover, in the definition of toroidal embeddings it is not necessary to make a base change to $\bar{k}$ if the embedding is strict \cite{Denef-ToroidalMorphisms}. In the definitions of toroidal embeddings and toroidal morphisms we can replace the completions by henselizations (see \cite{KKMS} page 195, and \cite{Denef-ToroidalMorphisms}). However, below we will not use the assertions made in the present remark.

\subsection{Toroidal actions}\label{sec-toroidal-actions}

An action of a finite group $G$ on a toroidal embedding $U\subset X$ defined over $\bar{k}$ is called a {\em toroidal action} at a closed point $x \in X$ if  there exist a local model $(V,v)$ at $x\in X$ and a group homomorphism $G_x\to T_{V, v}$ from the stabilizer $G_x$ of $x$ to the stabilizer $T_{V, v}$ of $v$ under the action on $V$ of the big torus $T_V\subset V$, such that the isomorphism
\[ \hat{\cO}_{X,x} \cong  \hat{\cO}_{V,v}\]
is compatible with the $G_x$-action, where $G_x$ acts on $\hat{\cO}_{V,v}$ through the homomorphism $G_x \to T_{V,v}$.
In this situation we say that the local model $(V,v)$ is \emph{compatible with} the $G$-action.
The action is {\em toroidal} if it is toroidal at every closed point.

Let $G$ act on $U\subset X$ toroidally, and assume the quotient $X/G$ exists.
Then the quotient $U/G\subset X/G$ is again a toroidal embedding (the local models are given by toric varieties $V/G_x$).

Let $f:(U_X\subset X) \to (U_B\subset B)$ be a toroidal morphism that is $G$-equivariant under toroidal actions of $G$ on both embeddings. We say that $f$ is {\em $G$-equivariantly toroidal} if for each closed point $x \in X$ there exist local models $(V,v)$ at $x$ and $(W,w)$ at $f(x)\in B$,
which are compatible with the morphism $f$ and the $G$-action on $X$ (hence also with the $G$-action on $B$).
Assuming again that the quotients $X/G$ and $B/G$ exist, such an $f$ induces a toroidal morphism of the quotient toroidal embeddings
\[ (U_X/G\subset X/G) \to (U_B/G\subset B/G).\]

For a toroidal embedding $U\subset X$ defined over $k$, an action of $G$ is called toroidal if the induced action on $U_{\bar{k}}\subset X_{\bar{k}}$ is toroidal.
The quotient of a toroidal embedding by a toroidal action is again a toroidal embedding. A $G$-equivariant morphism  $f:(U_X\subset X) \to (U_B\subset B)$ is called $G$-equivariantly toroidal if the base extension to $\bar{k}$ is $G$-equivariantly toroidal. Such a morphism induces a toroidal morphism of the quotient toroidal embeddings.

\subsection{Semistable families of curves.} \label{sec-semistable}

The reference here is \cite{deJong}.

A flat morphism $f:X\to B$ over the field $k$ is a {\em semistable family of curves} if every geometric fiber of $f$ is a complete reduced connected curve with at most ordinary double point singularities.

Consider  a semistable family of curves $f:X\to B$ over $\bar{k}$. If $x\in X$ is in the singular locus of $f$, then $X$ has a local equation at $x$:
\[  \hat{\cO}_{X,x} \cong \hat{\cO}_{B,f(x)}[[u,v]]/(uv-h)\]
for some $h\in \hat{\cO}_{B,f(x)}$. Note that $h=0$ defines the image of the singular locus of $f$ in $B$. It follows from this that if $U\subset B$ is a toroidal embedding and $f:X\to B$ is a semistable family of curves, smooth over $U$, then $f^{-1}(U)\subset X$ is also a toroidal embedding and the map $f$ is toroidal.

A similar statement holds for a semistable family of curves $f:X\to B$ defined over $k$: Suppose $U\subset B$ is a toroidal embedding and $f$ is smooth over $U$, then
\[ f: (f^{-1}(U)\subset X) \to (U\subset B)\]
is a toroidal morphism of toroidal embeddings.

\subsection{Resolution of singularities}

We will use the following two assertions
about resolution of singularities. They follow directly from any one of the canonical resolution of singularities algorithms \cite{Hironaka, EV, BM}, but are much weaker. Let $k$ be any field of characteristic zero.

\subsubsection{Assertion 1} \label{sec-resolution1}
Let $X$ be a variety over $k$, and $Z \subset X$ a proper closed subset. Then there exists a projective birational morphism $m_X: X' \rightarrow X$, defined over $k$, with $X'$ nonsingular, such that $Z' := m_X^{-1}(Z)$ is a strict normal crossings divisor on $X'$, and the restriction of $m_X$ to the complement of a suitable strict normal crossings divisor $Z'' \supset Z'$ is an open embedding.

\subsubsection{Assertion 2} \label{sec-resolution2}
Let $U \subset X$ be a toroidal embedding defined over $k$. Then there exists a projective birational morphism $m_X: X' \rightarrow X$, defined over $k$, with $X'$ nonsingular, such that $U' := m_X^{-1}(U) \subset X'$ is a strict toroidal embedding over $k$, and $m_X: (U' \subset X') \rightarrow (U \subset X)$ is toroidal.

\medskip
The first assertion follows directly from canonical resolution of singularities, but can
also be obtained by adapting the proof of weak resolution in \cite{AJ}, as we do in section \ref{sec-proof} below.

The second assertion is also a consequence of canonical resolution. Indeed, first apply canonical resolution to find a projective modification $m_1: X_1 \rightarrow X$, with $X_1$ nonsingular. Then $U_1 := m_1^{-1}(U) \subset X_1$ is a toroidal embedding and
$m_1: (U_1 \subset X_1) \rightarrow (U \subset X)$ is toroidal. To see this, we may assume that $X$ is toric, because canonical resolution commutes
with change of base field and with formal isomorphisms (see \cite{Kollar-Resolution}, Remark 3.56).
Moreover canonical resolution of a toric variety preserves the torus action, yielding an equivariant torus embedding as resolution space.
Next we apply canonical resolution to obtain an embedded resolution $m_2: X' \rightarrow X_1$ of $X_1 \setminus U_1$ in $X_1$. Then $U' := m_2^{-1}(U_1) \subset X'$ is not only a nonsingular toroidal embedding, but also strict, because $X' \setminus U'$ has strict normal crossings being the inverse image of $X_1 \setminus U_1$. One can avoid applying resolution in its full generality, replacing it by a combinatorial argument using logarithmic geometry, see \cite{Niziol-ToricSing}, Corollary 5.7. See also Remark \ref{remark-strict} below.

\section{Proof.} \label{sec-proof}
The purpose of this section is to prove Theorem \ref{th-toroidal-reduction}.

\subsection{Reduction to the projective case}


We first blow up $Z$ on $X$ and replace $Z$ by its inverse image; therefore we may assume $Z$ is the support of an effective Cartier divisor. This modification is projective in Grothendieck's sense; below we further modify $X$ by a quasi-projective variety, so the composite modification will be also projective in Hartshorne's sense \cite{HartshorneAG}.

Let us now reduce to the case where both $X$ and $B$ are quasi-projective varieties.
By Chow's lemma (\cite{EGA-II}, 5.6.1) there exist projective modifications $m_X:X'\to X$ and $B'\to B$ such that both $X'$ and $B'$ are quasi-projective varieties. Replacing $X'$ with the closure of the graph of the rational map $X'\das B'$, we may assume that $X'\to B'$ is a morphism, and that $m_X$ is still a projective morphism. Indeed the closure of the graph is contained in $X' \times_B B'$, and is hence projective over $X'$. Let $Z''$ be the union of $m_X^{-1}(Z)$ and the locus where $m_X$ is not an isomorphism.

Now we reduce to the {\em projective} case.
Choose projective closures $X'\subset \overline{X}$ and $B'\subset \overline{B}$, and again by the graph construction, we may assume that $\overline{X}\to \overline{B}$ is a morphism. Let $\overline{Z}=Z''\cup (\overline{X}\setmin X')$. Then the theorem for $\overline{Z}\subset \overline{X}\to \overline{B}$ implies it for $Z\subset X\to B$.
Indeed, given modifications $\overline{X}'\to \overline{X}$ and  $\overline{B}'\to \overline{B}$, we restrict these to the open sets $X'\subset   \overline{X}$ and $B'\subset   \overline{B}$. The inverse image of $Z$ in the modification of $X'$ is a divisor because $Z$ is the support of a Cartier divisor.
Thus, we may assume that $X$ and $B$ are projective varieties. Replacing $Z$ by a larger subset we may  assume at the same time that $Z$ is the support of an effective Cartier divisor on $X$.

\subsection{Structure of induction}\label{sec-inductionStruct} We assume that $X$ and $B$ are projective varieties, and proceed by induction on the relative dimension $\textrm{dim}\, X - \textrm{dim}\, B$ of $f:X\to B$. In the proof we will repeatedly replace $X$ and $B$ with suitable projective modifications, to which $f$ extends, until the requirements of the theorem are satisfied. This is certainly permitted if we replace $Z$ by a proper closed subset of the modification of $X$, which contains the inverse image of $Z$ and the locus where the modification is not an isomorphism. We will always (often without mentioning) replace $Z$ in this way, taking it large enough so that it is
the support of an effective Cartier divisor.

\begin{lem}(Abhyankar, cf. \cite{GM})\label{lem-abh} Let $X$ be a normal variety
and $f:X\rar B$ a finite
surjective  morphism onto a nonsingular variety, unramified outside a
divisor $D$ of normal crossings. Then $(B\setmin D\subset B)$ and $(X\setmin
f^{-1}(D)\subset X)$ are toroidal embeddings and $f$ is a toroidal
morphism. Moreover, if $f:X\rar B$ is Galois with Galois group $G$, then $G$ acts toroidally on $X$, and the stabilizer subgroups of $G$ are abelian. \qed
\end{lem}

\subsection{Relative dimension 0.} \label{sec-relDim0} Assume that the relative dimension of $f$
is zero. The proof in this case is a reduction to the Abhyankar's lemma.
\subsubsection{Constructing a Galois alteration}
We start the reduction by  constructing a projective alteration $\tilde{X}\to X$, such that
$\tilde{X}\to B$ is a Galois alteration. We assume that $X$ and $B$ are projective. Since $f$ is surjective, it is generically
finite. Let $L$ be a normalization of the function field $K(X)$ of $X$ over
the function field $K(B)$ of $B$. Then $L/K(B)$ is a finite Galois extension
with Galois group $G$. We choose a projective model $\tilde{X}$ of $L$ such
that $G$ acts on $\tilde{X}$.
For $\tilde{X}$ we can take for example the closure of the image of $X^0 \rightarrow \bar{X}^{|G|}: x \mapsto (g(x))_{g \in G}$, where $\bar{X}$ is a projective closure of an affine model $X^0$ of $L$ on which $G$ acts.
For each $g\in G$ we get a rational map $\phi_g:\tilde{X}\rar X$ corresponding to
the embedding $g: K(X) \hookrightarrow L$ obtained by restricting $g$ to $K(X)$. We let
\[ \overline{\Gamma}\subset \tilde{X}\times \prod_{g\in G} X \]
be the closure of the graph of $\prod_{g\in G} \phi_g$. Then the group $G$
acts on $\overline{\Gamma}$ and projection to one of the factors $X$ gives us
a morphism $\overline{\Gamma}\rar X$. So, we may replace $\tilde{X}$ by
$\overline{\Gamma}$ and assume that the rational map $\tilde{X}\rar X$ is a
morphism.
\subsubsection{Reduction to a finite morphism}
Consider the quotient variety $\tilde{X}/G$. Since $G$ fixes the field
$K(B)$, we have a birational morphism  $p:\tilde{X}/G\rar B$, hence a rational
map $p^{-1}\circ f: X\rar \tilde{X}/G$. Let $X_0\subset X\times_B
\tilde{X}/G$ be the closure of the graph of this rational map, and note that
the projection $\tilde{X}\rar\tilde{X}/G$ factors through $X_0$:
\[
\begin{array}{rcl}
    &     &  \tilde{X} \\
    &  \swarrow & \dar \\
X_0 & \rar & X \\
f_0\dar & & \dar f \\
\tilde{X}/G& \stackrel{p}\rar & B
 \end{array}
\]
In the diagram above the horizontal maps are modifications.
Since $\tilde{X}\rar \tilde{X}/G$ is finite, so is $f_0: X_0 \rar
\tilde{X}/G$. Thus, we have modified $f:X\to B$ to a finite morphism $f_0: X_0\to \tilde{X}/G$. It suffices now to prove the theorem for $f$ replaced by $f_0$, with $Z$ replaced by a proper closed subset $Z_0$ of $X_0$ satisfying the requirements of section \ref{sec-inductionStruct}.

Alternatively, the reduction to a finite morphism can also be proved by applying the Flattening Theorem of Raynaud and Gruson \cite{Raynaud-Gruson}, and the fact that any proper flat generically finite morphism is finite.

\subsubsection{Applying  Abhyankar's lemma}
Let $D\subset \tilde{X}/G$ be the branch locus of $f_0$. We let
$B'\rar \tilde{X}/G$ be a resolution of singularities (see \ref{sec-resolution1}) such that the inverse
image of $D\cup f_0(Z_0)$ is contained in a strict normal crossings divisor $D'$ on $B'$, with $B' \rightarrow \tilde{X}/G$ an isomorphism outside $D'$. Denote by $X'$ the irreducible component of the normalization of $X_0\times_{\tilde{X}/G} B'$ that dominates $X_0$. Then the projection $f': X'\rar B'$ is a finite
morphism unramified outside the normal crossings divisor $D'$. By Abhyankar's lemma such a morphism is toroidal with respect to the toroidal embeddings $\; U_{B'}:=B' \setminus D' \subset B'\;$ and $\; U_{X'}:=f'^{-1}(U_{B'}) \subset X'$.

Note that, by construction, $U_{B'}\subset B'$ is a nonsingular strict toroidal embedding. Applying resolution of singularities (see \ref{sec-resolution2}) to $X'$ and its divisor $X'\setmin U_{X'}$, we may assume that  $U_{X'}\subset X'$ is also a nonsingular strict toroidal embedding. The morphism $f': (U_{X'}\subset X') \to(U_{B'}\subset B')$ and $Z'\subset X'$ satisfy the statements of the theorem. This finishes the proof of the theorem in case $f$ has relative dimension 0.   \qed

\medskip
Assume now that we have proved the theorem for morphisms of relative dimension $n-1$, and consider the case that $f$ has relative dimension $n$, with $X$ and $Y$ projective varieties, $f$ surjective, and $Z$ the support of an effective Cartier divisor.

\subsection{Preliminary reduction steps}

The idea of the proof in case of relative dimension $n$ is to factor the morphism $f:X\to B$ as a composition $X\to P\to B$, where $X\to P$ has relative dimension $1$ and $P\to B$ has relative dimension $n-1$. We then apply the induction assumption to the morphism $P\to B$, after having replaced $X\to P$ by a semistable family of curves (section~\ref{sec-semistable}), using semistable reduction. In order to apply the semistable reduction theorem \cite{deJong-families}, we need the map $X\to P$ to have geometrically irreducible generic fiber. Let us construct such a factorization.

\subsubsection{Normalizing}
First, we may replace $X$ with its normalization, therefore we can
assume $X$ is normal, replacing $Z$ as explained in section \ref{sec-inductionStruct}.
Let $\eta\in B$ be the generic point of $B$. Similarly, we can assume that $B$ is normal.

\subsubsection{Using Bertini's theorem} By the
projectivity assumption we have $X\subset \bP^N_B$ for some
$N$.
Let $L\subset \bP^N_\eta$ be a general enough $(N-n)$-plane, so that $L\cap X$ is finite and contained in the nonsingular locus of $X$, and such that no line in $L$ is tangent to $X$. The set of such $L$ contains a nonempty open subset $U_1$ of the Grassmannian $\GG(N-n, \bP^N_\eta)$ of $(N-n)$-planes in $\bP^N_\eta$.
Let $\bP^N_\eta \das \bP^{n-1}_\eta$ be the projection from  $L\subset \bP^N_\eta$.  This gives a rational map $X_\eta \das \bP^{n-1}_\eta$ that is not defined at the finite set of points $L\cap X$. Blowing up these points gives a projective morphism $\tilde{X}_\eta\to   \bP^{n-1}_\eta$ with fibres  $M\cap X$ for all $(N-n+1)$-planes $L\subset M$.
Indeed the blowup $\tilde{X}_\eta$ is the closure of the graph of
$X_\eta \das \bP^{n-1}_\eta$. Note that $\tilde{X}_\eta$ is normal because $X$ is normal and the center of the blowup is a finite set of nonsingular points. Since $X$ is normal, a general enough $(N-n+1)$-plane in $\bP^N_\eta$ is disjoint from the singular locus of $X$. Thus by Bertini's Theorem (see e.g. \cite{HartshorneAG}, Chapter III, Corollary 10.9 and Remark 10.9.2) there exists a nonempty open subset $U_2$ of the Grassmannian $\GG(N-n+1, \bP^N_\eta)$ of $(N-n+1)$-planes in $\bP^N_\eta$, such that the scheme-theoretic intersection $M\cap X$ is nonsingular for each $M \in U_2$. Let $\Gamma$ be the closed subset of
$\GG(N-n, \bP^N_\eta) \times \GG(N-n+1, \bP^N_\eta)$ consisting of the pairs $(\alpha, \beta)$ with $\alpha \subset \beta$. Note that $\Gamma$ is irreducible, since it is an image of an open subset of an affine space. Hence, the image of the projection
$(U_1 \times U_2)\cap \Gamma \rightarrow U_1$ is dense in $U_1$, because the projections of $\Gamma$ on the two Grassmannians are surjective. We conclude that
there exists a nonempty open set $U_1' \subset U_1$ of planes in the Grassmannian $\GG(N-n, \bP^N_\eta)$, such that the generic fibre of  the morphism $\tilde{X}_\eta\to   \bP^{n-1}_\eta$ is smooth, whenever $L \in U_1'$. Because the field $k$ is infinite, the $k$-valued points are dense in the Grassmannian, hence we may choose the plane $L$ to be defined over $k$.

\subsubsection{Using Stein factorization} The rational map $X_\eta \das \bP^{n-1}_\eta$ gives a rational map $X\das \bP^{n-1}_B$, defined over $k$. Let us replace $X$ with the normalization of the closure of the graph of this map, so we may assume we have a morphism   $X\to \bP^{n-1}_B$, with $X$ normal. The generic fibre of this morphism is smooth (it is the same as the generic fibre of $\tilde{X}_\eta\to   \bP^{n-1}_\eta$).
Let $X\stackrel{g}{\rar}P\rar\bP^{n-1}_B$ be
the Stein factorization, where $g:X\rar P$ is a projective morphism of relative dimension $1$ with geometrically connected fibers, and the second morphism is finite (see \cite{EGA-III}, 4.3.1 and 4.3.4).
Then the generic fibre of $g$ is geometrically irreducible. Since $X$ is normal, also $P$ is normal.

We are now ready to apply the semistable reduction theorem to the morphism $g$.

\begin{definition} Let $\alpha: X_1\rar X$ be an alteration with $X$ a variety over $P$, and $Z\subset X$ an
irreducible divisor dominating $P$. The {\em strict altered transform} $Z_1\subset X_1$ of $Z$ is
the closure of $\alpha^{-1}(\eta)$ in $X_1$, where $\eta$ is the generic
point of $Z$. The strict altered transform of an arbitrary divisor is the union of the strict altered transforms of its components that dominate $P$.
\end{definition}

\subsection{Semistable reduction of a family of curves}

By \cite{deJong-families}, Theorem 2.4, items (i)-(iv) and (vii)(b), there exists a commutative diagram of morphisms of normal projective varieties

$$\begin{array}{lcl} X_1 & \stackrel{\alpha}{\rar} & X \\ \dar g_1 & &
                     \dar g \\ P_1 & \stackrel{a}{\rar} & P \\ & &
                     \dar \\ & & B \end{array}
$$
 and a finite group $G\subset \Aut_PP_1$, with the following
 properties:

\begin{enumerate}
\item The morphism $a:P_1\to P$ is a Galois alteration with Galois
group $G$ (i.e. $P_1/G\rar P$ is birational).
\item The action of $G$ lifts to $\Aut_XX_1$, so that $g_1$ is $G$-equivariant, and $\alpha:X_1\to X$ is
a Galois alteration with Galois group $G$.
\item There are disjoint sections $\sigma_i:P_1\rar X_1$, $i=1, \dots,\kappa,$ such that
the strict altered transform $Z_1\subset X_1$ of $Z$ is the union of
their images and $G$ permutes the sections $\sigma_i$.
\item The morphism $g_1: X_1 \to P_1$ is a semistable family of curves with smooth generic fibre, and
$\sigma_i(P_1)$ is disjoint from $\operatorname{Sing}g_1$ for each $i$.
\end{enumerate}

Note that the image of $\alpha^{-1}(Z)\setminus Z_1$ in $X$ lies over a proper closed subset of $P$, because $X$ has relative dimension 1 over $P$. The same holds for the locus in $X_1/G$ where the modification $X_1/G \rightarrow X$, induced by $\alpha$, is not an isomorphism. Indeed $X$ is normal, hence any rational map from $X$ to a complete variety is regular outside a subset of codimension $\geq 2$. Thus we can find an effective Cartier divisor on $X_1/G$ whose support $Z'$ contains $\alpha^{-1}(Z)/G$, such that $X_1/G \rightarrow X$ is an isomorphism outside $Z'$, and
$Z' \setminus (Z_1/G)$ lies over a proper closed subset of $P_1/G$.

We may replace $X$, $P$, $Z$ by $X_1/G$, $P_1/G$, and
$Z'$, cf. the discussion in section \ref{sec-inductionStruct}.
Then $X_1/G = X$, $P_1/G = P$, and $\alpha^{-1}(Z)\setminus Z_1$ lies over a proper closed subset of $P_1$. Finally, observe that the singular locus of $g_1: X_1 \rightarrow P_1$ lies over a proper closed subset of $P_1$, since $g_1$ is flat (because semistable) with smooth generic fibre (by (4)).

\subsection{Using the inductive hypothesis}
Let $\Delta\subset P$ be the union of the loci over which $P_1$ or
$X_1$ are not smooth, and the closure of the image of $\alpha^{-1}(Z)\setminus Z_1$ in $P$. Note that $\Delta$ is a proper closed subset of $P$.
We apply the inductive assumption to
$\Delta\subset P \to B$, and obtain a diagram
$$\begin{array}{lclcl} U_{P'} & \hookrightarrow & P' &\stackrel{m}{\to} &P \\
\dar                      &                      & \dar & & \dar\\
 U_{B'} & \hookrightarrow & B' &\to &B
                     \end{array}
$$
such that $P' \to P$ and $B' \to B$ are
projective modifications, the left square is a toroidal morphism of nonsingular strict toroidal embeddings,
$m^{-1}\Delta$ is a divisor of strict normal crossings contained in
$P' \setmin U_{P'}$, and $m$ is an isomorphism on $U_{P'}$.

We may again replace $P, B$ by $P', B'$, writing $U_P, U_B$ instead of $U_{P'}, U_{B'}$, and further we may replace
$X,X_1, P_1, \sigma_i$ by the normalizations of their pullbacks to $P'$, and $Z$ by the union of its inverse image and the inverse image of $P' \setmin U_{P'}$. With the \emph{pullback} to $P'$ of a variety over $P$, we mean here the irreducible component of the base change to $P'$ that dominates the given variety.
After these replacements the properties (1), (2), (3), (4) and the equalities $X_1/G = X$, $P_1/G = P$, are still true.
Indeed, the base change of the old $X_1$ to the new $P_1$ equals the new $X_1$, since it is integral (by flatness) and normal (by semistability).
Moreover, both $\alpha^{-1}(Z)\setminus Z_1$ and the singular locus of $g_1$ lie over $P \setminus U_P$.
Now, $P$ and $B$ are nonsingular and $P \setminus U_P$ is a strict normal crossings divisor on $P$.
Moreover $P\to B$ is toroidal, and $P_1\to P$ is
unramified over $U_P$.

By Abhyankar's Lemma~\ref{lem-abh}, since $P_1$ is normal, it inherits
a toroidal structure given by $U_{P_1} = a^{-1}(U_P)$ as well, so that $P_1\to P$
is a finite toroidal morphism.

To summarize, in addition to properties (1)-(4) above, and the equalities $X_1 / G = X , $ $P_1/G = P$, we also have that the morphisms $a: (U_{P_1}\subset P_1) \to  (U_P\subset P)$ and  $(U_P\subset P) \to (U_B\subset B)$ are toroidal. Moreover the toroidal embeddings $U_P \subset P$ and $U_B \subset B$ are strict and nonsingular. The embedding $U_{X_1}\subset X_1$, where $U_{X_1} = g_1^{-1}U_{P_1}\setmin \cup_i \sigma_i(P_1)$, is a toroidal embedding and the morphism $g_1:  (U_{X_1}\subset X_1) \to  (U_{P_1}\subset P_1)$ is a $G$-equivariant toroidal morphism, by section \ref{sec-semistable} and (4). Note also that the divisor $\alpha^{-1}(Z)$ lies in $X_1\setmin U_{X_1}$.

\subsection{Torifying the group action}\label{sec-torification}

 By Abhyankar's Lemma \ref{lem-abh}, $G$ acts toroidally on $(U_{P_1} \subset P_1)$ and its stabilizers are abelian. If $G$ acts toroidally on $(U_{X_1} \subset X_1)$ and if the morphism $g_1: X_1 \rightarrow P_1$ is $G$-equivariantly toroidal, then the induced morphism $X_1/G \rightarrow P_1/G = P$ is toroidal, cf. section \ref{sec-toroidal-actions}. Moreover, by resolution  of singularities (see \ref{sec-resolution2}) we find a nonsingular strict toroidal embedding $U_{X'} \subset X'$ and a toroidal modification $X' \rightarrow X_1/G$.
We then obtain a toroidal morphism $X' \rightarrow X_1/G \rightarrow P \rightarrow B =: B'$, as required by Theorem \ref{th-toroidal-reduction}.
However, in general, $G$ does not act toroidally on $X_1$.

We follow section 1.4 of \cite{AJ} to construct a suitable modification of $X_1$ on which $G$ acts toroidally. In \cite{AJ} this modification is obtained by two blowups, each followed by normalization. However, there one works over an algebraically closed field $\bar{k}$, thus we need to verify that the ideals blown up are actually defined over $k$, so that the modification is also defined over $k$.
We recall the construction of the ideals to be blown up and explain why they are defined over k. For the convenience of the reader we also recall why these constructions yield a toroidal action, although this is all done in \cite{AJ}.

\subsubsection{Blowing up the singular locus} A first situation where $G$ does not act toroidally on $X_1$ happens when an element of $G_x$ exchanges two components of a fiber $g_{1,\bar{k}}: X_{1,\bar{k}} \rightarrow P_{1,\bar{k}}$ passing through a point $x \in X_{1,\bar{k}}$. This problem is solved in \cite{AJ} by blowing up the singular scheme $S$ of the morphism $g_{1,\bar{k}}$, hence separating all nodes. Note that $S$ is the subscheme of $X_{1,\bar{k}}$ defined by the first Fitting ideal sheaf of $g_{1,\bar{k}}$. This ideal sheaf is obtained from the first Fitting ideal sheaf of $g_1$ by base change. Thus $S$ is defined over $k$. Let $X_2$ be the blowup of $X_1$ along $S$, and $X_2^{\textrm{nor}}$ the normalization of $X_2$. The action of $G$ on $X_1$ lifts to an action of $G$ on $X_2$ and $X_2^{\textrm{nor}}$. Let $U_{X_2}$ be the inverse image of $U_{X_1}$ in $X_2$. Note that $U_{X_2}$ is nonsingular because the morphism $g_1$ is smooth on $U_{X_1}$. We identify the inverse image of $U_{X_1}$ in $X_2^{\textrm{nor}}$ with $U_{X_2}$.

\subsubsection{Local description} \label{sec-localDescription} First we recall why $U_{X_2} \subset X_2^{\textrm{nor}}$ is a toroidal embedding. Let $x$ be a closed point of $X_{2,\bar{k}}^{\textrm{nor}}$.
Choose a local model $(V,v)$ of $P_{1,\bar{k}}$ at the image of $x$ in $P_{1,\bar{k}}$, compatible with the $G$-action as in section \ref{sec-toroidal-actions}, and let $A$ be the completion of the local ring of the toric variety $V$ at $v$. As can be seen from the local description (section \ref{sec-semistable}) of the semistable family $X_{1,\bar{k}}$ over $P_{1,\bar{k}}$, there are only 3 possible cases for the completion of $X_{2,\bar{k}}$ at the image of $x$ in $X_{2,\bar{k}}$, namely one of the following formal spectra:
$$
\textrm{Spf} \, A[[y,z]]/(zy^2-h),\; \textrm{Spf} \, A[[y,z]]/(y^2-h),\, \textrm{or}\;\, \textrm{Spf} \, A[[z]],
$$
where in the first two cases $h \in A$ is a monomial (i.e. a character of the big torus of $V$), and $h(x)=0$.
The third case holds if and only if the image of $x$ in $X_{1,\bar{k}}$ is not in the singular locus of $g_{1,\bar{k}}$. Only in this case it is possible that $x$ belongs to the inverse image $\Gamma$ of $\cup_i \sigma_i(P_1)$ in $X_{2,\bar{k}}$, and then we may assume that $z=0$ is a local equation for $\Gamma$ at $x$.
The first formal spectrum and the third one are completions of appropriate (not necessarily normal) equivariant torus embeddings. The same holds for each component of the formal spectrum of the normalization of $A[[y]]/(y^2-h)$. Hence $U_{X_2} \subset X_2^{\textrm{nor}}$ is a toroidal embedding. Note also that the divisor $X_{2,\bar{k}}^{\textrm{nor}} \setminus U_{{X_2},\bar{k}}$ is $G$-strict.

\subsubsection{Analyzing the group action} In the first case, the ideal generated by $y$, as well as the ideal generated by $z$, is invariant under the action of $G_x$. Hence multiplying $y$ and $z$ by suitable units with residue 1, we may assume that $G_x$ acts on $y$ and $z$ by characters of $G_x$. Indeed, replace $y$ by
$|G_x|^{-1} \sum_{\sigma \in G_x} (y/\sigma(y))(x)\, \sigma(y)$.
Thus the action of $G$ on $X_{2,\bar{k}}^{\textrm{nor}}$ is toroidal at $x$.

In the third case, if $x \in \Gamma$ then a similar argument as in the first case shows that the action of $G$ on $X_{2,\bar{k}}^{\textrm{nor}}$ is toroidal at $x$.

However, in the second and third case, if $x \notin \Gamma$ then
the action on $U_{{X_2},\bar{k}} \subset X_{2,\bar{k}}^{\textrm{nor}}$ is in general not toroidal at $x$, but because $G_x$ is abelian, we can choose the local formal parameter $z$ such that $G_x$ acts on it by a character (indeed consider the representation of $G_x$ on the  vector space over $\bar{k}$ generated by the $G$-orbit of $z$). When this character is nontrivial, the action is not toroidal at $x$. Indeed, the divisor locally defined at $x$ by $z=0$ is not contained in the toroidal divisor $X_{2,\bar{k}}^{\textrm{nor}} \setminus U_{{X_2},\bar{k}}$. Moreover these locally defined divisors, as $x$ varies, might not come from a globally defined divisor (because these are not defined in a canonical way).

\subsubsection{Pre-toroidal actions} At any rate, the action of $G$ on $X_{2,\bar{k}}^{\textrm{nor}}$ is pre-toroidal in the sense of Definition 1.4 of \cite{AJ}. A faithful action of a finite group $G$ on a toroidal embedding $U \subset X$ over $\bar{k}$ is called \emph{pre-toroidal} if the divisor $X \setminus U$ is $G$-strict and if for any
closed point $x$ on $X$ where the action is not toroidal we have the following. There exists an isomorphism $\epsilon$, compatible with the $G_x$-action and the toroidal structure, from the completion of $X$ at $x$ to the completion of
$X_0 \times \textrm{Spec}\, \bar{k}[z]$
at $(x_0,0)$, where $U_0 \subset X_0$ is a toroidal embedding with a faithful toroidal $G_x$-action, $x_0$ is a closed point of $X_0$ fixed by $G_x$, the toroidal structure on
$X_0 \times \textrm{Spec}\, \bar{k}[z]$
is given by
$U_0 \times \textrm{Spec}\, \bar{k}[z] \subset X_0 \times \textrm{Spec}\, \bar{k}[z]$ and the action of $G_x$ on
$X_0 \times \textrm{Spec}\, \bar{k}[z]$
comes from its action on $X_0$ and its action on z by a nontrivial character $\psi_x$ of $G_x$. Note that the character $\psi_x$ only depends on $x$ and not on $\epsilon$. Assume now that the $G$-action on $U \subset X$ is pre-toroidal.

\subsubsection{Blowing up the torific ideal} \label{sec-torificBlowup} In \cite{AJ} (Theorem 1.7 and the proof of Proposition 1.8) it is proved that there exists a canonically defined $G$-equivariant ideal sheaf $\cI$ on $X$, called the \emph{torific ideal sheaf}, having the following properties. At each closed point of the support of $\cI$ the $G$-action is not toroidal. For each closed point $x$ of $X$, where the $G$-action is not toroidal, the completion of $\cI$ at $x$ is generated by the elements of $\widehat{\cO}_{X,x}$ on which $G_x$ acts by the character $\psi_x$. And finally, if we denote by $\widetilde{X}$ the normalization of the blowup of $X$ along $\cI$, and by $\widetilde{U}$ the inverse image of $U$ in $\widetilde{X}$, then $\widetilde{U} \subset \widetilde{X}$ is a toroidal embedding on which $G$ acts toroidally.

The proof of this last property follows directly from the following local description at $x$ of the blowup, where we may assume that $X_0$ is an affine toric variety, $U_0$ is the big torus of $X_0$, $G_x$ is a subgroup of $U_0$, $G_x$ acts on $X_0$ through $U_0$,
$X = X_0 \times \textrm{Spec}\, \bar{k}[z]$,
$U = U_0 \times \textrm{Spec}\, \bar{k}[z]$,
and $\epsilon$ is the identity.
Locally at $x$ the ideal sheaf $\cI$ is generated by $z$ and monomials $t_1, \cdots, t_m$ in the coordinate ring $R$ of $X_0$. Indeed, at least one monomial of $R$ is contained in $\cI_x$, because $G_x$ is a subgroup of the big torus of $X_0$. Above a small neighborhood of $x$, the blowup of $X$ along $\cI$ is covered by the charts
$$
\textrm{Spec} \, R[z,t_1/z, \cdots, t_m/z],  \; \textrm{Spec} \, R[t_i,t_1/t_i,\cdots, t_m/t_i][z/t_i],
$$
for $i=1, \cdots, m$. These are torus embeddings of $U_0 \times \textrm{Spec}\, \bar{k}[z, z^{-1}]$, hence their normalizations are toric. Moreover the embedding $\widetilde{U} \subset \widetilde{X}$ is toroidal, and $G_x$ acts toroidally on it, at any closed point of $\widetilde{X}$ above $x$, because on the first chart the locus of $z=0$ is contained in the inverse image of $X \setminus U$, and on the other charts the action of $G_x$ on $z/t_i$ is trivial. The above argument also shows that the support of $\cI$ is disjoint from $U$, because $\cI_x$ contains a monomial in $R$. Thus the blowup is an isomorphism above $U$.

If the toroidal embedding $U \subset X$ and the $G$-action are defined over $k$, then the torific ideal sheaf $\cI$ is also defined over $k$, because it is stable under the action of the Galois group of $\bar{k}$ over $k$. Indeed this is a direct consequence of the above description of the completions of $\cI$. Hence $\widetilde{X}$ is also defined over $k$.

\subsubsection{Conclusion of the proof} We now apply this to $X_2^{\textrm{nor}}$. Let $X_3$ be the blowup of $X_2^{\textrm{nor}}$ along the torific ideal sheaf, $X_3^{\textrm{nor}}$ the normalization of $X_3$, and $U_{X_3}$ the inverse image of $U_{X_2}$ in $X_3$. The blowup morphism $X_3 \rightarrow X_2^{\textrm{nor}}$ is an isomorphism above $U_{X_2}$. Hence $U_{X_3}$ is nonsingular, and we identify it with the inverse image of $U_{X_2}$ in $X_3^{\textrm{nor}}$. Note that $G$ acts toroidally on the toroidal embedding $U_{X_3} \subset X_3^{\textrm{nor}}$.
Moreover the morphism $X_3^{\textrm{nor}} \rightarrow X_2^{\textrm{nor}} \rightarrow X_1$ is an isomorphism above $U_{X_1}$.
By the argument in the beginning of section \ref{sec-torification} (which uses \ref{sec-resolution2}) with $X_1$ replaced by $X_3^{\textrm{nor}}$,
we now see that in order to prove Theorem \ref{th-toroidal-reduction} it suffices to show that the composition of the morphisms
$X_3^{\textrm{nor}} \rightarrow X_2^{\textrm{nor}} \rightarrow X_1 \rightarrow P_1$ is $G$-equivariantly toroidal. But this is a straightforward consequence of the local description of $X_{2,\bar{k}}$ given in \ref{sec-localDescription}, and the local description of the blowup of $X$ along $\cI$ given in \ref{sec-torificBlowup}, with $X = X_{2,\bar{k}}^{\textrm{nor}}$. This finishes the proof of Theorem  \ref{th-toroidal-reduction}.  \qed

\subsubsection{Remark} \label{remark-strict}
It follows from section 1.4.3 of \cite{AJ} that the divisor $D_3 := X_3^{\textrm{nor}} \setminus U_{X_3}$ is $G$-strict, when $k$ is algebraically closed. Hence, in that case, the toroidal embedding $U_{X_3}/G \rightarrow X_3^{\textrm{nor}}/G$ is strict, and one can apply the toroidal resolution of \cite{KKMS} (chapter II \S2) to desingularize the quotient $X_3^{\textrm{nor}}/G$, instead of applying assertion \ref{sec-resolution2} as we do above. This is the road followed in \cite{AJ,Abramovich-Karu}. Actually, if we would use the fact that $D_{3,\bar{k}}$ is $G$-strict, then in order to prove Theorem \ref{th-toroidal-reduction}, we only need assertion \ref{sec-resolution2} in the special case that the given toroidal embedding becomes strict after base change to $\bar{k}$. In that case the method in chapter II of \cite{KKMS} can be used to desingularize, and afterwards strictness can be achieved as explained in section 2.4 of \cite{deJong}.

\bibliographystyle{plain}
\bibliography{ADK}


\end{document}